\documentclass{article}
\usepackage{latexsym,amsmath,geompsfi}
\newcommand{\be}{\begin{equation}}
\newcommand{\ee}{\end{equation}}
\newcommand{\ra}{\rightarrow}

\newcommand{\R}{{\bf R}}
\newcommand{\Z}{{\bf Z}}
\newcommand{\Atop}[2]{\genfrac{}{}{0pt}{}{#1}{#2}}
\newcommand{\0}{^{\phantom0}}
\newcommand{\makefig}[3]{
        \begin{figure}[htbp]
        \refstepcounter{figure}
        \label{#2}
        \begin{center}
                ~#3~\\
                \medskip
                {\sf Figure \thefigure:  #1}
        \end{center}
        \end{figure}
}
\begin{document}

\begin{Large}
\centerline{On some points-and-lines problems and configurations}
\end{Large}

\vspace*{1ex}

\centerline{Noam D. Elkies}

\vspace*{2ex}

{\bf Abstract.}  We apply an old method for constructing
points-and-lines configurations in the plane
to study some recent questions in incidence geometry.

\begin{quote}
  What are known as ``Points and Lines'' puzzles
  are found very interesting by many people.
  The most familiar example, here given,
  to plant nine trees so that they shall form
  ten straight rows with three trees in every row,
  is attributed to Sir Isaac Newton, but the earliest collection
  of such puzzles is, I believe, in a rare little book that I possess
  --- published in 1821 --- {\em Rational Amusement for Winter Evenings},
  by John Jackson.  The author gives ten examples of
  ``Trees planted in Rows.''

  These tree-planting puzzles have always been a matter of great perplexity.
  They are real ``puzzles,'' in the truest sense of the word, because
  nobody has yet succeeded in finding a direct and certain way
  of solving them.
  They demand the exercise of sagacity, ingenuity, and patience,
  and what we call ``luck'' is also sometimes of service.

  --- H.E.~Dudeney, \em{Amusements in Mathematics}\/ (1917)
  \cite{Dudeney2}, page~56
\end{quote}

{\bf Introduction.}
Almost a century after Dudeney wrote these paragraphs,
problems in incidence geometry continue to perplex
both recreational and professional mathematicians, and
the prospect of a uniform ``direct and certain way of solving them''
remains remote.  Even for natural asymptotic questions,
a wide gap often separates the best upper and lower bounds known.
In this paper we construct some explicit point-and-line configurations 
that yield new lower bounds for two specific questions of this kind.
Question~1, suggested by the recreational literature, asks:
How many lines can meet $n^2$ points in the plane
in at least $n$ points each?  Question~2 arises in the research
literature~\cite{Brass}:  If on each of $N$\/ horizontal lines
we choose (at most) $N$\/ points, how many additional lines
can contain $N$\/ of these $N^2$ points?
It turns out that an arrangement of $16$ points in $15$ lines of~$4$
(Figure~\ref{fig:n=4} below), which has been known at least since 1908,
naturally generalizes to configurations that
not only give lower bounds for Question~1
but also improve on the previous records for Question~2.
We also find a variation of this construction that yields
a partial answer to Question~1 and a further improvement for the cases
$N=12m=12,24,36,\ldots$ and $N=12m-1=11,23,35,\ldots$ of Question~2.
By the construction in~\cite{Brass}, the new results for Question~2
yield, for each $N \geq 5$, improved lower bounds on the exponent
in the asymptotic ``orchard-planting'' problem
with \hbox{$N$-point} lines.
Each of these arrangements exploits dihedral symmetry:
the lines include all axes of symmetry,
and every point lies on one of the axes
and at least one pair of lines symmetrical with respect to this axis.
This approach is at least a century old
(we give specific citations later), but might still produce
further new examples and results for modern incidence geometry.

The rest of this paper is organized as follows.
We first give some general background
on this kind of points-and-lines problem.
We then introduce Question~1, on plane arrangements
of $n^2$ points with many \hbox{$n$-point} lines,
and show the best configurations previously known.
Next we present Brass's problem as Question~2,
and observe that some of the configurations
already known for Question~1 also answer Question~2.
We proceed to modify the known constructions
to obtain further improvements for both Questions.
Finally we reconsider the symmetry of our configurations,
which can be even greater than it appears.  Most notably,
the obvious fivefold dihedral symmetry of Figure~\ref{fig:n=4}
extends to an action of the icosahedral group $A_5$
by projective linear transformations.
This action, and an analogous action
of the octahedral group $S_4$ on the real projective plane,
leads us to further points-and-lines configurations
related with the finite projective planes of orders $\leq 5$.
We expand the customary concluding Acknowledgements, to explain
how we became aware of Question~2 and its connection with Question~1
even though such problems are quite far from our usual research work.

{\bf Definitions of $T_k$ and $t_k$, and of $T_k^{(r)}$ and $t_k^{(r)}$;
the exponents $\tau_k$.}
For a finite set~$S$\/ of points in the plane,
let $t_k(S)$ $(k=2,3,4,\ldots)$
be the number of lines meeting~$S$\/ in exactly $k$ points,
and $T_k(S) = \sum_{k' \geq k} t_{k'}(S)$
the number of lines meeting~$S$\/ in at least $k$ points.
For a positive integer~$n$ let
\be
t_k(n) := \max_{|S|=n} t_k(S),
\quad\
T_k(n) := \max_{|S|=n} T_k(S),
\label{eq:Tk}
\ee
so $t_k(n)$ (or $T_k(n)$) is the largest number of lines
that can contain exactly (or at least) $k$ points
out of some configuration of $n$ points in the plane.
Clearly $t_k(n) \leq T_k(n)$.
For $r>k$ we also let
\be
\def\scs{\scriptstyle}
t_k^{(r)}(n) := \max_{\Atop{\scs |S|=n}{\scs T_r\0\!(S)=0}} t_k(S),
\quad\
T_k^{(r)}(n) := \max_{\Atop{\scs |S|=n}{\scs T_r\0\!(S)=0}} T_k(S),
\label{eq:Tkr}
\ee
restricting $S$\/ to point sets for which
no line contains $r$ or more points.
For instance, the condition that no line contain more than~$k$\/ points
(common in ``orchard-planting'' problems) corresponds to $r=k+1$,
and clearly in that case $t_k^{(r)}(n) = T_k^{(r)}(n)$.
See for instance \cite[p.315~ff.]{BMP},
where $t_k^{(k+1)}(n)$ is called $t_k^{\rm orchard}(n)$.

A key question concerns the asymptotic behavior
of $t_k^{(r)}(n)$ and $T_k^{(r)}(n)$
as $n \ra \infty$ for fixed~$k,r$.
The question is trivial for $k=2$:
clearly $t_2^{(3)}(n) = \binom{n}{2} = T_2^{(r)}(n)$ for all $r>2$.
In general, for all $k,r,n$ we have an elementary upper bound
$T_k^{(r)}(n) \leq \binom{n}{2} \big/ \binom{k}{2}$.
For $k=3$ this gives $T_3^{(r)}(n) \leq n^2/6 - O(n)$,
which is known to be asymptotically sharp:
certain configurations of torsion points on cubic curves
even give $t_3^{(4)}(n) = n^2/6 - O(n)$
(see for instance \cite{Burr,BGS}).
For $k \geq 4$, Erd\H{o}s proposed long ago the conjecture that
$t_k^{(k+1)}(n) = o(n^2)$
(this is ``Conjecture~12'' of \hbox{\cite[p.317]{BMP})};
more generally one might guess that $T_k^{(r)}(n) = o(n^2)$
for any fixed $k,r$ with $4 \leq k < r$.
[Note that the corresponding conjecture for $T_k(n)$ or even $t_k(n)$
is false, for instance because a $k\times m$\/ lattice array
has at least $m^2/(k-1)$ lines of exactly $k$\/ points
(and even this is not optimal, see~\cite{Palasti});
this is why we fix some finite upper bound~$r$
on the number of points in any line.]
But it is not known that $t_k^{(k+1)} = o(n^2)$ for any $k\geq 4$,
even though the best lower bounds on $T_k^{(r)}(n)$ are
$C_{k,r} n^{\tau_k}$ with
\be
2 = \tau_3 > \tau_4 \geq \tau_5 \geq \tau_6 \geq \cdots \ra 1.
\label{eq:tau}
\ee
Our results include improvements on these $\tau_k$ for each $k \geq 5$
(though to be sure we are still nowhere near settling
Erd\H{o}s's conjecture).
See Theorem~1, stated near the end of this paper.

{\bf Sets of $n^2$ points in the plane with many \hbox{$n$-point} lines.}
Anyone who has seen a magic square knows that
$t_n^{(n+1)}(n^2) \geq 2n+2$ for all $n>1$:
a square array of $n^2$ points in the plane
forms $2n+2$ lines of~$n$, namely
the $n$ horizontal lines, $n$ vertical lines, and $2$ diagonals.
For $n=2$ this is clearly optimal because each of the six pairs
of points has a two-point line through it.
But for each $n>2$ one can get more than $2n+2$ lines.
A famous configuration, known at least since the beginning
of the twentieth century \cite[p.175]{Dudeney1},
shows that for $n=4$ one may get as many as $15$ lines of~$4$
by using a double pentagram instead of a square.
See Figure~\ref{fig:n=4}.  (The closed and open circles indicate
points on $3$ and $5$ lines respectively; more about this later.)

\makefig{$16$ points, $15$ lines of~$4$}
{fig:n=4}{\psfig{figure=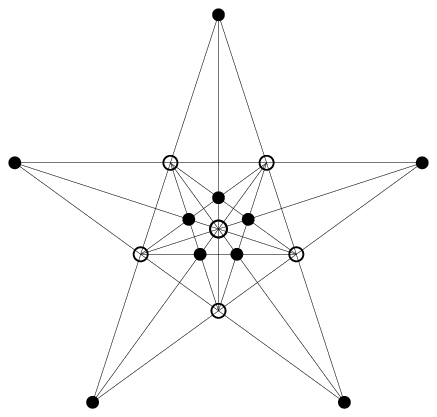,width=144pt} \\
}

This construction readily generalizes to all even $n>2$:
replace the two nested pentagrams by two nested \hbox{$(n+1)$-point} stars,
each formed from the longest diagonals of a regular \hbox{$(n+1)$-gon},
to obtain a configuration with \hbox{$(n+1)$-fold} dihedral symmetry
consisting of $n^2$ points lying on $3n+3$ lines of $n$ points each.
Figure~\ref{fig:n=6} shows the case $n=6$ of this construction.

\makefig{$36$ points, $21$ lines of~$6$}
{fig:n=6}{\psfig{figure=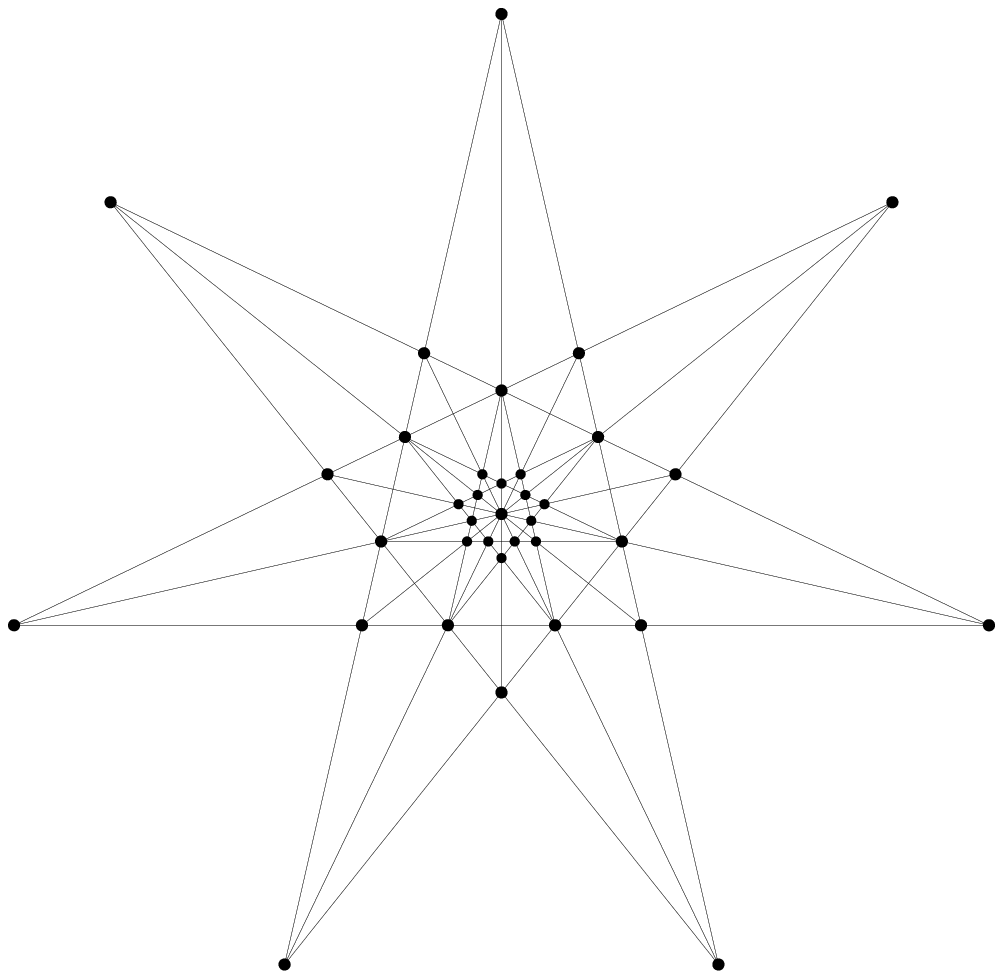,width=252pt} \\
}

\makefig{$9$ points, $10$ lines of~$3$}
{fig:n=3}{\psfig{figure=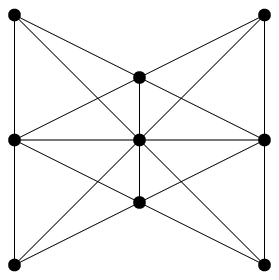,width=72pt} \\
}

This suggests several questions, which we first raised
in the interview~\cite[p.228]{Vakil}:

{\bf Question 1a}\/: {\em Is this configuration optimal?}

That is, is $3n+3$ the maximal number of lines that can meet $n^2$ points
in the plane in (at least) $n$ points each?
Using the notation of~(\ref{eq:Tk}), we are asking:
is $t_n(n^2) = T_n(n^2) = 3n+3$ for $n=4,6,8,\ldots$?
This might be known for $n=4$,
but is almost certainly open for every even $n\geq 6$.

{\bf Question 1b}\/: {\em What happens for odd~$n$?}

For $n=3$ it has long been known that the maximum is~$10$,
though over the complex numbers the famous configuration of
nine flex points of a smooth cubic has~$12$ lines of three
(it is probably mere coincidence that this is also
the value of $3n+3$ for $n=3$),
which attains the upper bound $\binom92\big/\binom32$ exactly:
the line through {\em every} pair of points goes through
a third point of the configuration.\footnote{
  Note too that the flexes of a smooth cubic in the plane
  are also its \hbox{$3$-torsion} points.
  Over the real numbers, we already noted the use of torsion points
  on such curves in estimating $t_3^{(4)}(n)$.
  For more on points-and-lines arrangements
  in the complex plane and beyond, see \cite{EPS,Kelly}.
  }
The \hbox{$10$-line} configuration, mentioned by Dudeney
in the passage quoted earlier from \cite[p.56]{Dudeney2},
is obtained from the $3 \times 3$ square array
by moving an opposite pair of edge points halfway towards the center
(Figure~\ref{fig:n=3});
we later return to this configuration as well.\footnote{
  Burr begins his article \cite{Burr} by quoting the puzzle asking
  for this configuration from the same source
  ({\em Rational Amusement for Winter Evenings} (1821) by John Jackson),
  where it is given as a verse:
\begin{tabbing}
  ab\=cde\=\kill
  \> Yo\=ur aid I want, nine trees to plant\\
  \> \>  In rows just half a score;\\
  \> And let there be in each row three.\\
  \> \>  Solve this: I ask no more.
\end{tabbing}
  }

Some twenty years ago we constructed --- with some ``luck''$\!$,
as Dudeney might say --- a sporadic arrangement of
$25$ points with $18$ lines of~five (Figure~\ref{fig:n=5},
also shown in~\cite[p.228]{Vakil}).
The points on each edge of the triangle bisect and trisect the edge.
Thus $t_n(n^2) \geq 3n+3$ also for $n=5$.
We construct a different such configuration
later, from Figure~\ref{fig:n=5+}.
We do not know whether $18$ lines is maximal,
nor whether either \hbox{$18$-line} configuration was known earlier.

\makefig{$25$ points, $18$ lines of~$5$}
{fig:n=5}{\psfig{figure=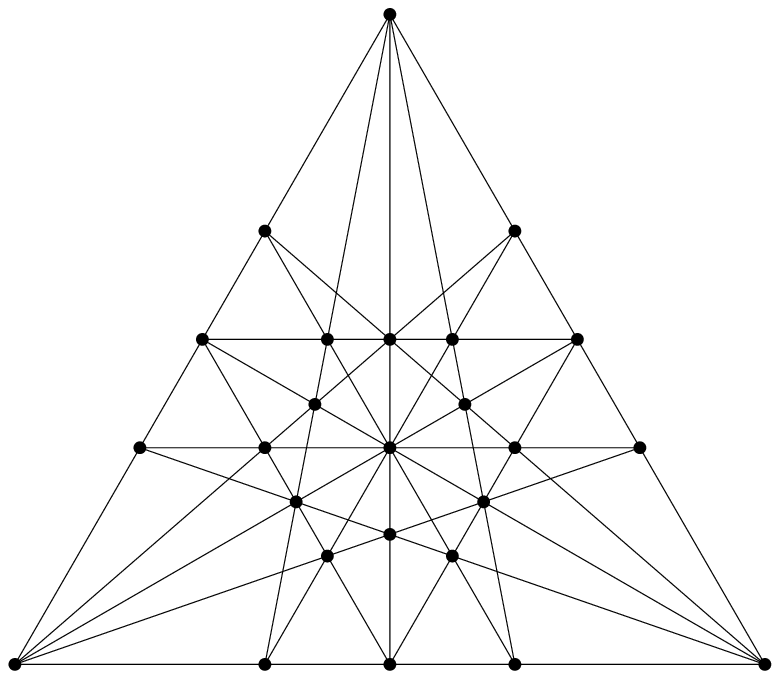,width=185pt} \\
}

For odd $n\geq 7$, one can at least see quickly that
the $2n+2$ lines of the square configuration are {\em not}\/ optimal.
We can already get $2n+2$ lines using only $(n-1)^2+3$ points:
$(n-1)^2$ in a square array, one point in the center of the square
(which has not been used yet because $n$ is odd),
and two points at infinity where the line at infinity
meets the coordinate axes.  Then we can use the remaining $2n-4$ points
to form another $4\sqrt{n}-O(1)$ lines by putting them on diagonals
that contain $n-2$, $n-3$, $n-4$, \ldots~points in the array.
At the end, if points at infinity are deemed undesirable one may apply
a projective transformation to put all $n^2$ points
in the finite plane.\footnote{
  Dudeney used much the same trick in his second solution
  \cite[p.190]{Dudeney2} to the puzzle of placing $21$ points
  in $12$ lines of~$5$: the configuration is projectively
  equivalent to the $16$ points
  $(i,j)$ ($1 \leq i,j \leq 4$) in the $(x,y)$ plane,
  together with the points $(3/2, 5/2)$ and $(5/2, 3/2)$
  and the three points at infinity contained in
  the four lines $x=i$, the four lines $y=j$,
  and the three lines $x=y$, \, $x-y=\pm 1$.
  The twelfth line is then $x+y=4$.
  The use of projections in this context to bring points at infinity
  to the finite plane is noted explicitly in \cite[p.105]{BallCox}.
  }
We shall show that $3n+1$ lines can always be attained,
even under the ``orchard'' constraint
that no line contain more than~$n$ points;
that is, $t_n^{(n+1)}(n^2) \geq 3n+1$.
We shall also show that for $n = 12m-1 = 11$,~$23$,~$35$,~$47$,~etc.,
there are configurations of $n^2$ points in the plane
with $3n+4$ lines each of which passes through at least $n$
of the $n^2$ points.  But these configurations necessarily contain
some lines of $n+1$ points, so we obtain $T_n(n^2) \geq 3n+4$
for these values of~$n$ but not $t_n(n^2) \geq 3n+4$.

{\bf Parallel lines with many Brass transversals.}
P.~Brass asks \cite{Brass}:

{\bf Question 2}\/:
{\em
Can there be $N$\/ parallel lines~$l_i$ in the plane,
and $M>N+4$ lines $\lambda_j$ not parallel to the~$l_i$,
such that for each~$i$ we have
\be
\# \Bigl( \bigcup_{j=1}^M l_i \cap \lambda_j \Bigr)  \leq N
\label{eq:brass_ineq}
\ee
(that is, there are at most~$N$\/ points on~$l_i$
through which some $\lambda_j$ passes)?
}

We shall call such $\lambda_j$ a collection of ``Brass transversals''
to the $l_i$.  More generally, one may of course ask,
for any $N$\/ and~$N'$, for the maximal number of lines~$\lambda_j$
whose union intersects each of $N$\/ parallel lines~$l_i$
in at most $N'$ points.
But the case $N=N'$ is of particular interest
because Brass~\cite{Brass} gives an explicit recursive construction
showing that a collection of $M$\/ Brass transversals yields
$t_N^{(N+1)}(n) = \Omega(n^{\log_N \! M})$ as $n \ra \infty$.

Question~2 specifies $M>N+4$ because $M=N+4$ can be attained
for each $N\geq 3$.  Let $l_i$ be the line $x=i$ for $i<N$,
and the line at infinity for $i=N$\/;
let $\lambda_j$ be the line $y=j$ for $j \leq N$\/;
and let the remaining four transversals be the lines
$y=x$, $y=x+1$, $x+y=N$, and $x+y=N+1$.
These $N+4$ lines meet $l_N$ in $3$ points,
and $l_i$ in~$N$\/ points for each $i<N$.
For $N=3$ this configuration is easily seen to be unique
up to projective transformations.
Figure~\ref{fig:n=3} shows it in another guise,
with $7=3+4$ Brass transversals to the three vertical lines;
projecting one of these lines to infinity yields the case $N=3$
of the construction described earlier in this paragraph.
The resulting bound $t_3^{(4)}(n) = \Omega(n^{\log_3 7})$
is not interesting, because we already know that
$t_3^{(4)}(n)$ is asymptotic to $n^2/6$.
But in \cite[p.317]{BMP} we find that
for $N \in [5,17]$ the lower bound with exponent $\log_N\0(N+4)$
is the best exponent known, and for $N \geq 18$ it can be used with
a different recursive construction due to Gr\"unbaum~\cite{Grunbaum}
to obtain the record exponent $1 + (1/(N-\gamma))$
with $\gamma \doteq 3.59$.

We improve this to
\be
\log_N\0 2N = 1 + \frac{\log 2}{\log N}
\label{eq:better_tau}
\ee
for each $N=5,7,9,\ldots$,
using our configurations from Question~1 with $n=N-1$ and $M=2N$.
Project the center of our \hbox{$n^2$-point} configuration
to infinity;
let the $l_i$ be the $N=n+1$ lines through this point at infinity,
and let the $\lambda_j$ be the remaining $M=2N=2n+2$ lines.
Then $\#(\cup_j l_i \cap \lambda_j) = N$\/ for each~$i$,
and the bound $t_N^{(N+1)}(n) = \Omega(n^\tau)$
with $\tau = \log_N\0 2N$\/ follows by~\cite{Brass}.

We cannot quite do this for $N=6$ using our sporadic
\hbox{$25$-point} configuration in Figure~\ref{fig:n=5},
because the six lines through the center are not equivalent.
When $l_i$ is one of the three axes of symmetry of the triangle,
the $\lambda_j$ meet $l_i$ in only four points;
but for the other three~$l_i$ (those parallel to the triangle's sides),
there are seven points of intersection.
Still, this configuration may be of use for Brass's construction
because the inequality~(\ref{eq:brass_ineq}) remains true on average,
even with a strict inequality:
one might have expected eight points of intersection
for $l_i$ in the second group, but the two new points coincide
because two of the $\lambda_j$ are parallel to~$l_i$
and thus meet $l_i$ in the same point at infinity
(which is not one of the $25$ points of our configuration).

{\bf Further refinements.}
For $N=3$ the configuration that attains $N+4=7$ Brass transversals
is unique, and can be displayed symmetrically as shown
on the left side of Figure~\ref{fig:2N+1}
by projecting one of the transversals to infinity.
This again suggests a generalization to arbitrary odd~$N$\/:
let $l_i$ be the line through the origin making angle $(i/N)\pi$
with the horizontal; and let $\lambda_j$ be the $N$\/ pairs of lines
parallel to the~$l_i$ at unit distance, together with
the line at infinity, for a total of $M=2N+1$ transversals.
Taking the indices of the $l_i$ modulo~$N$, we see that
for each $i' \bmod N$\/ the transversals parallel to $l_{i\pm i'}$
meet $l_i$ at the point(s) $1/\sin((i'/N)\pi)$ units from the origin.
This gives $N$\/ points of intersection for each line,
and all the points with $i'=0$ are on the line at infinity,
which accounts for the \hbox{$(2N+1)$-st} transversal.
The right side of Figure~\ref{fig:2N+1}
shows the $N=5$ case of this construction.
Again we conclude by projecting the origin to infinity
to obtain parallel lines~$l_i$.
For each $N=5,7,9,\ldots$,
this gives us an even better value $\log_N\0(2N+1)$
for the exponent $\tau_N\0$ of~(\ref{eq:tau}).
Moreover, the set of $N^2$ points $l_i \cap \lambda_j$
meets the $3N+1$ lines $l_i, \lambda_j$ in $N$\/ points each,
and meets no line in more than~$N$\/ points
because the set is contained in the $N$\/ lines~$l_i$.
Therefore $t_N^{(N+1)}(N^2) \geq 3N+1$.
This gives a new lower bound on $t_N^{(N+1)}(N^2)$
for each odd $N \geq 7$.  (We exclude $N=5$,
because then $3N+1=16$, but Figure~\ref{fig:n=5} already attains~$18$.)

\makefig{$n$ lines, $2n+1$ Brass transversals
  including the line at infinity ($n=3,5$)}
{fig:2N+1}{\psfig{figure=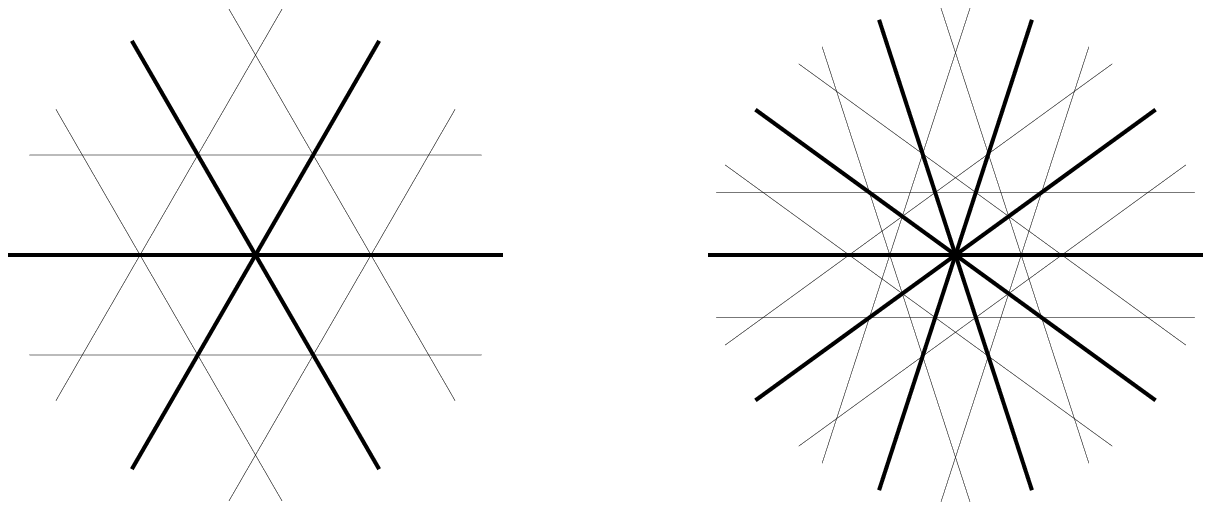,width=312pt}\\
}

This construction fails when $N$\/ is even,
because then the points at unit distance from the origin on~$l_i$
each lie on just one transversal (with $i'=N/2$).
But we still achieve $M=2N$\/ by discarding the line at infinity
and rotating the other lines $\lambda_j$ by an angle $\pi/2N$
about the origin.  This improves on $N+4$ for all even $N\geq 6$.
(Figure~\ref{fig:2N} shows the case $N=6$.)  We therefore attain
$\tau_N\0 = \log_N\0 2N$\/ for all even~$N$, and have thus improved
the exponent $\tau_N\0$ for all integers $N \geq 5$.

Our configuration with a double \hbox{$N$\/-point} star
also required that $N$\/ be odd, for a different reason:
for even~$N$, the longest diagonals of a regular \hbox{$N$\/-gon}
that do not go through its center intersect each other in only
$N-2$ points.  But for large~$N$\/ the double-star construction
has some flexibility that we can sometimes exploit
to improve the configuration and allow some even~$N$\/ as well.
Namely, we may match any of one star's rings of $N$\/ intersection points
with any ring at a different position on the other star.
This can be done when the ratio between stars' circumradii is
$$
\rho_N\0(i,j) := \sin \frac{i\pi}{N} \Bigl/ \sin \frac{j\pi}{N} \Bigr.
$$
for some distinct positive integers $i,j < N/2$,
regardless of the parity of~$N$.
[So far, as in Figure~\ref{fig:n=6} (with $N=7$),
we have always used $(i,j)=(1,(N-1)/2)$.]
If $\rho_N\0(i,j) = \rho_N\0(i',j')$
for another pair $(i',j')$ of integers in $(0,N/2)$,
then the resulting double-star configuration
has the same number of incidences with $N$\/ fewer points.
We may find such $i,i',j,j'$ when $6|N$\/ and $N \geq 12$,
using the identity
$$
\sin \theta \, \sin (\frac{\pi}{2} - \theta )
= \sin \theta \, \cos \theta
= \frac12 \sin 2\theta
= \sin \frac{\pi}{6}  \,\sin 2\theta.
$$

\makefig{$6$ lines, $12$ Brass transversals}
{fig:2N}{\psfig{figure=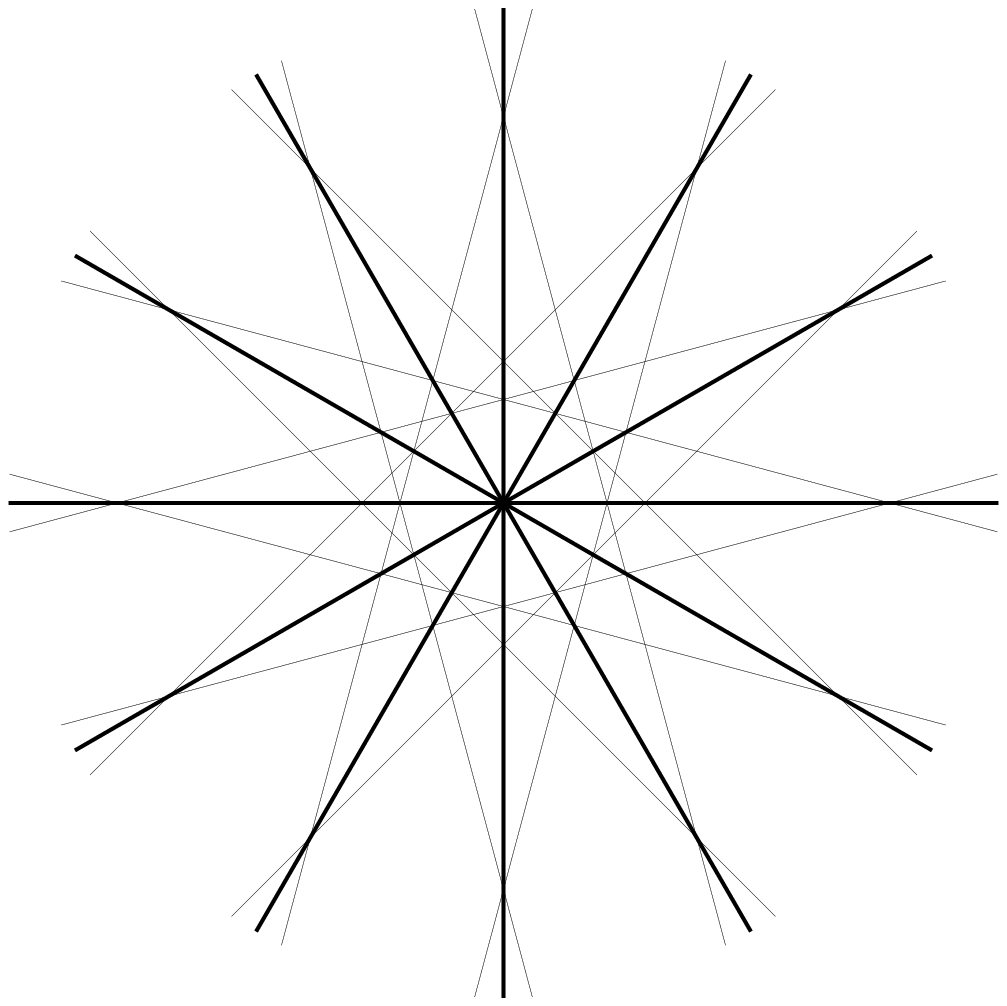,width=288pt}\\
}

[That these are in fact the only solutions is a special case
(and much easier than the full result) of \cite[Thm.~4]{PR};
the authors of~\cite{PR} report that the same theorem
had already been obtained by Bol~\cite{Bol}.
Unfortunately it is not possible to have a third pair $(i'',j'')$.]
This gives $(i',j,j') = (2i, N/6, (N/2)-1)$.
Moreover, when $N=12m$, we may choose $i,i',j,j'$ so that
$i$ and~$j'$ are odd while $i'$ and~$j$ are even, for instance
$(i,i',j,j')=(1,2,2m,6m-1)$.
Figure~\ref{fig:n=11} shows this when $N=12$.
Projecting the center to infinity then yields $N$\/ parallel lines
and $2N+1$ Brass transversals
(including the projection of the line at infinity, as before),
with only $N-1$ intersection points on each parallel line.
We have thus obtained yet another improvement
for the cases $N=12m$ and $N=12m-1$ of Question~2.

We can also use this configuration to partly answer Question~1b,
as follows.  Each of the lines through the center has a pair of points
each of which lies on just one of the transversals.
(These $N$\/ pairs of points
are marked by closed circles in Figure~\ref{fig:n=11}.)
There are $2^N$ sets of $N$\/ points
containing one point from each of these pairs;
choosing one of these sets and removing it leaves $(N-1)^2$ points
with $2N$\/ lines of $N-1$ points and $N+1$ lines of~$N$.
We have thus shown that $T_n(n^2) \geq 3n+4$ for
$n = N-1 = 12m-1$, as promised  earlier.

\makefig{$(N,M) = (11,25)$ or $(12,25)$;
also, $121$ points, $37$ lines of at least~$11$ (see text)}
{fig:n=11}{\psfig{figure=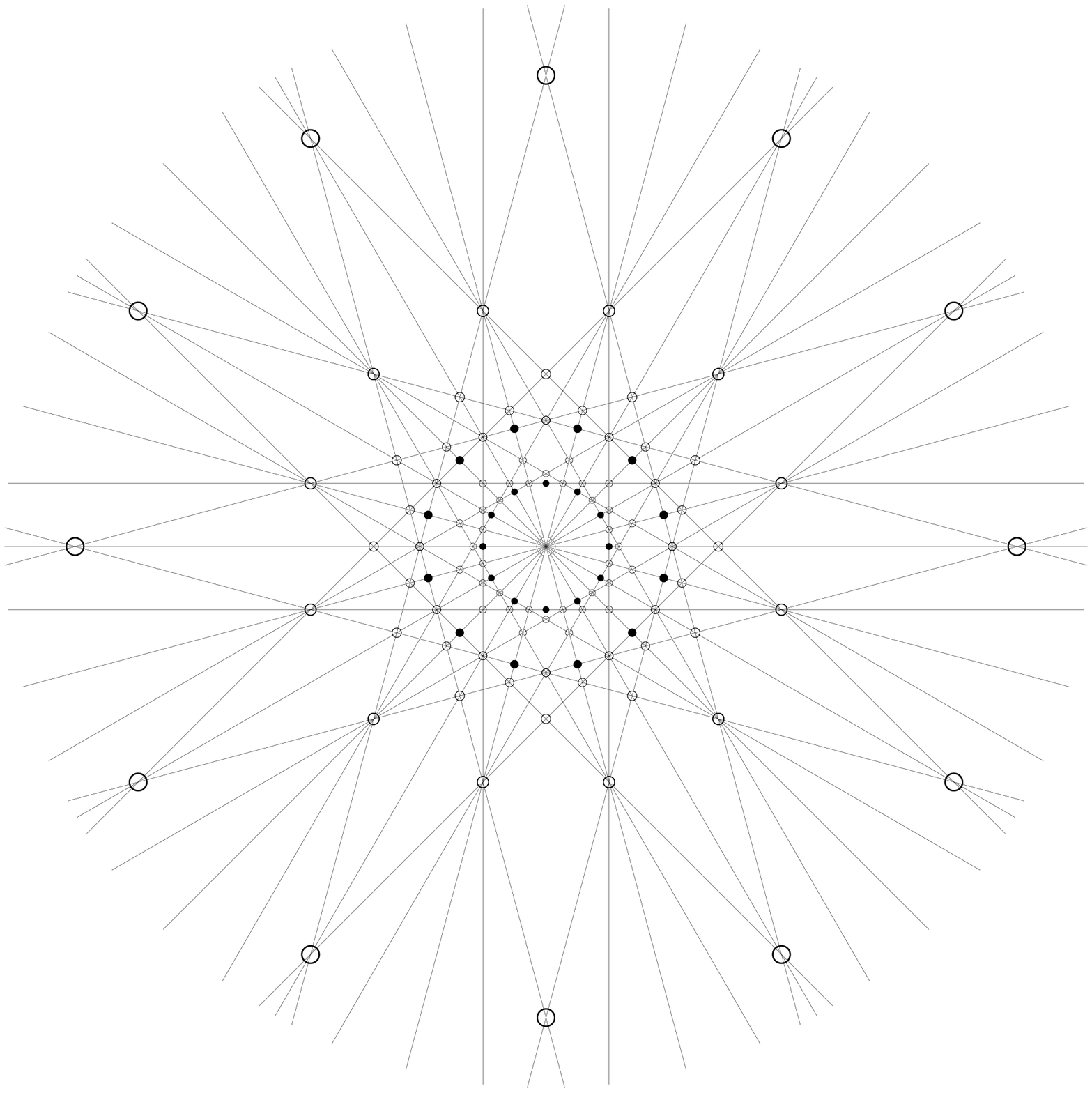,width=\textwidth}\\
}

Returning to Question~2, we collect all our results
and use them in Brass's recursive construction~\cite{Brass},
obtaining:

{\bf Theorem~1.}  {\em
For $N \geq 5$, let
\be
M = M(N) =
\begin{cases}
2N + 3, & \text{ if $N \equiv -1 \bmod 12$; } \\
2N + 1, & \text{ if $N \equiv 0, 1, 3, 5, 7$, or~$9 \bmod 12$; } \\
2N    , & \text{ if $N \equiv 2, 4, 6, 8$, or~$10 \bmod 12$. } \\
\end{cases}
\label{eq:M}
\ee
Then for each $N\geq 5$ we have $t_N^{(N+1)}(n) = \Omega(n^{\tau_N})$
where $\tau_N\0 = \log_N M$.
}

A numerical table of these new $\tau_N\0$ for $5 \leq N \leq 30$ follows:

\addtolength{\arraycolsep}{-2pt}

\vspace*{1ex}

\centerline{
$
\begin{array}{c|ccccccccccccc}
N & 5 & 6 & 7 & 8 & 9 & 10 & 11 & 12 & 13 & 14 & 15 & 16 & 17
\\ \hline
\tau_N\0 & 1.489 & 1.386 & 1.391 & 1.333 & 1.340 & 1.301 &
  1.342 & 1.295 & 1.284 & 1.262 & 1.268 & 1.250 & 1.254
\end{array}
$
}

\vspace*{2ex}

\centerline{
$
\begin{array}{c|ccccccccccccc}
N & 18 & 19 & 20 & 21 & 22 & 23 & 24 & 25 & 26 & 27 & 28 & 29 & 30
\\ \hline
\tau_N\0 & 1.239 & 1.244 & 1.231 & 1.235 & 1.224 & 1.241 &
1.224 & 1.221 & 1.212 & 1.215 & 1.208 & 1.210 & 1.203
\end{array}
$
}

\vspace*{1ex}

\centerline{Table 1}

These values of $\tau_N\0$, like the ones previously known,
approach~$1$ as $N \ra \infty$, but much more slowly:
$\tau_N\0 - 1 \approx \log 2 / \log N$,
while the previous results had $\tau_N\0 - 1 \approx 1/N$.
Unlike those previous $\tau_N\0$, the values in Table~1
are quite far from the monotonic descent described in~(\ref{eq:tau}).
For instance, our lower bound on $t_6^{(r)}(n)$ (any $r>7$)
uses configurations with many \hbox{$7$-point} lines,
and for $N=8,9,10$ and $r>11$ our lower bound on $t_N^{(r)}(n)$
uses configurations with many \hbox{$11$-point} lines!
Evidently the asymptotic behavior of $t_k^{(r)}(n)$ remains
``a matter of great perplexity''$\!$,
as Dudeney described it almost $90$ years ago.
Can one improve on Theorem~1 by showing that
$\tau_N\0 \geq \tau_{N'}\0$ when $N \leq N'$?
Can one exploit the extra line $l_{N+1}$
in our configuration for the case $N=12m-1$ of Question~2
to obtain a further asymptotic improvement?
Can the configurations arising from the identity
$\rho_N\0(i,2i) = \rho_N\0(N/6,(N/2)-1)$
be exploited also in the cases $N=18,30,42,\ldots$ when $N/6$ is odd?

One can attempt similar constructions with three or more nested stars,
or only one.  The only such variation we have found
that bears on the questions that motivated us here is a triple pentagram.
Adding to the old \hbox{$16$-point} configuration
of Figure~\ref{fig:n=4} a third star, and also each of the five points
where the line at infinity meets parallel sides of the three stars,
we obtain $26$ points spanning $21$ lines of~$15$.  
See Figure~\ref{fig:n=5+}.
The open circles mark the $10$ points
each of which is contained in only three of the $21$ lines;
removing any one of these leaves $25$ points in $18$ lines of~$5$,
in a configuration distinct from Figure~\ref{fig:n=5}.

\makefig{$26$ points ($5$ at infinity), $21$ lines of~$5$}
{fig:n=5+}{\psfig{figure=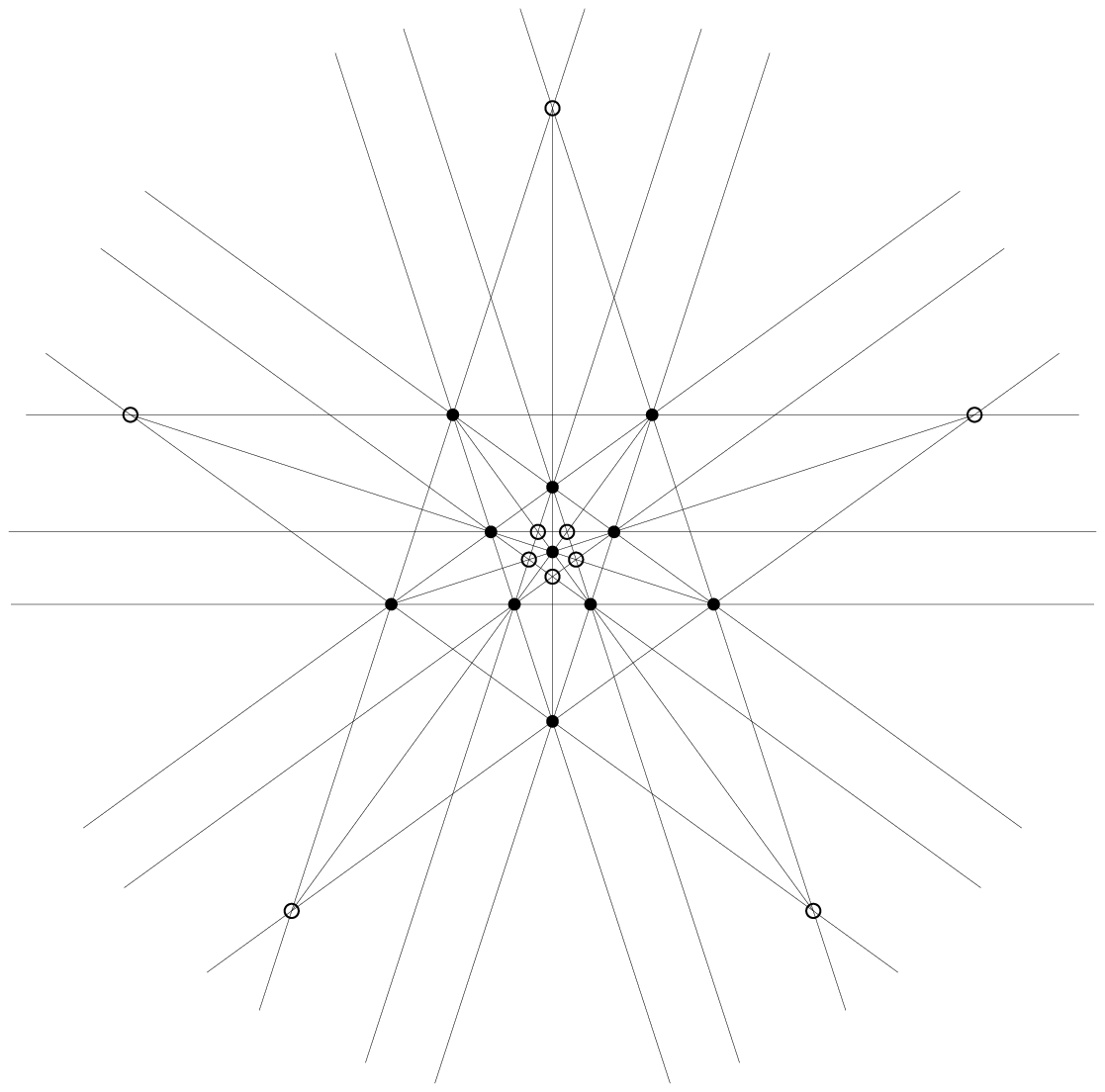,width=320pt} \\
}

\makefig{$13$ points ($3$ at infinity), $13$ lines}
{fig:n=3b}{\psfig{figure=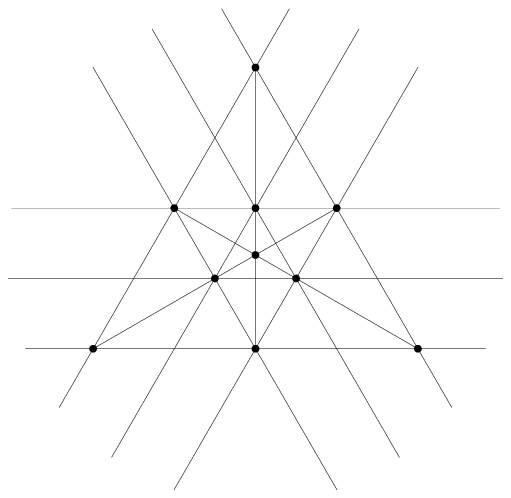,width=144pt} \\
}

{\bf More about Figure~\ref{fig:n=4} and symmetries.}
We saw that the solution of the puzzle ``to plant nine trees so that
they shall form ten straight rows with three trees in every row''
is more symmetrical than it appears from its usual presentation
in Figure~\ref{fig:n=3}:
this presentation has only $4$ symmetries,
but the projection shown on the left side of Figure~\ref{fig:2N+1}
exhibits the \hbox{$12$-element} group of symmetries of the regular hexagon.
Likewise, our initial configuration of $16$ points in $15$ lines of~$4$
(Figure~\ref{fig:n=4}) turns out to be even more symmetrical
than it looks: its group of projective symmetries
is the alternating group $A_5$,
acting transitively on the $15$ lines
and dividing the points into orbits of size $6$ and~$10$.
(The six-point orbit consists of
the central point and the five points of the middle ring,
each of which lies on $5$ four-point lines;
these are the points drawn as open circles in Figure~\ref{fig:n=4}.)
To see this, let $A_5$ act on the $2 \cdot 6$ vertices
of a regular icosahedron in~$\R^3$,
and map those vertices to $6$ points in~${\bf P}^2$
while preserving a fivefold symmetry of the icosahedron.
The other $10$ points are the pairs of face centers,
and the $15$ lines are dual to the pairs of edge centers.\footnote{
  We noted this online at \cite{Elkies}.  This page links
  to a picture of the images of the vertices, face centers,
  and edge centers, and of their dual lines;
  it also mentions Question~1 and the configurations for
  the cases $n=4,6,8,\ldots$ and $n=5$.
  }
Let $P_6$, $P_{10}$, and $P_{15}$ be the
$6$-, $10$-, and \hbox{$15$-point} orbits of points
under this action of~$A_5$,
and $L_6$, $L_{10}$, $L_{15}$ the corresponding orbits of lines.
Then for $i,j\in\{6,10,15\}$ there exists a point in $P_i$
contained in some line of $L_j$ if and only if $i=15$ or $j=15$,
in which case there are $30$ such points.
Figure~\ref{fig:n=4} shows $60$ of these incidences.
If we instead consider the $21$ points and $21$ lines of
$P_6 \cup P_{15}$ and $L_6 \cup L_{15}$, we find $90$ incidences.
These are contained among the $105$ incidences
in the finite projective plane of order~$4$;
the $15$ missing incidences are between each point of~$P_{15}$
and its dual line.  Likewise the $31$ points and $31$ lines of
$P_6 \cup P_{10} \cup P_{15}$ and $L_6 \cup L_{10} \cup L_{15}$
show $150$ of the $186$ incidences of the projective plane of order~$5$.

A similar configuration arises from the regular cube or octahedron,
with symmetry group $S_4$,
again larger than can be shown in any plane projection.
The vertices, faces and edges of a regular octahedron yield
$3+4+6$ points and as many lines, shown in Figure~\ref{fig:n=3b}.
The $48$ incidences
are among the $52$ in the finite projective plane of order~$3$,
lacking only the incidences between each face point and its dual line.
To explain this, note that the points are the images of
the $2 \cdot 13$ nonzero points $(x_1,x_2,x_3) \in \Z^3$
with each $|x_i|\leq 1$, and likewise the lines are
\hbox{$x_1 y_1 + x_2 y_2 + x_3 y_3 = 0$}
with each $y_i\in\{-1,0,1\}$.
These remain distinct when reduced mod~$3$,
and the only new incidences mod~$3$ are the four with
$$
 (x_1 : x_2 : x_3) = (y_1 : y_2 : y_3) = (\pm 1 : \pm 1 : \pm 1) .
$$
We can similarly relate the configurations of $21$ or $31$
points and lines of the previous paragraph with the corresponding
finite projective planes,
by recognizing them as points and lines with small coordinates
in $\Z[\varphi]$ where $\varphi = (\sqrt{5}+1)/2$,
and then reducing these coordinates modulo the prime ideal
$2\Z[\varphi]$ or $\sqrt5 \Z[\varphi]$ respectively.

{\bf Account and acknowledgements.}
Last year I traveled to Calgary for the Workshop in Discrete Geometry
in honor of the 50th birthday of K\'aroly Bezdek, and found my way
to the lecture room just in time for the problem session.
I intended to present an open ``tree-planting'' problem
in incidence geometry (Question~1) that I had wondered about
for some time.  The first few cases lead to appealing configurations;
I had no better reason than pure curiosity for asking the question
in general, but this meeting seemed a natural venue to raise
the problem, and a reasonable one to hope for new information.
That incidence geometry was an appropriate topic was confirmed
when Peter Brass, who was among the first to present a problem
at this session, asked a question of a similar flavor (Question~2),
though his interest in it was more than recreational:
a positive answer would yield an asymptotic improvement
to a construction in his recent paper~\cite{Brass}.
I thought that one of the ``appealing configurations''
I was about to show (Figure~\ref{fig:n=4}) might work,
and after some hurried scribbling verified that
projecting its center point to infinity
answers the first odd instance ($N=5$) of Brass's question.
Later experimentation showed that the natural generalization
of this configuration (as in Figure~\ref{fig:n=6} for $n=6$)
yields such an answer for all odd $N>3$, and afterwards led to the
further refinements described in the Introduction and illustrated
in Figures~\ref{fig:2N+1}, \ref{fig:2N}, and~\ref{fig:n=11}.

I thank the organizers of the Calgary Workshop in Discrete Geometry,
for inviting me to participate in the workshop;
Peter Brass, for extended e-mail correspondence
on these problems, including references to his paper~\cite{Brass}
and the relevant sections of~\cite{BMP}; and
the referee, for directing me to references \cite{Ismailescu, Palasti}
and suggesting a rearrangement of the exposition.
This paper is based on research supported in part
by NSF grant DMS-0501029.

{\small
Dept.\ of Mathematics\\
Harvard University\\
Cambridge, MA 02138, USA\\ \\
{\tt elkies@math.harvard.edu}
}

\begin{thebibliography}{9}
\bibitem{BallCox} Ball, W.W.Rouse; Coxeter, H.S.M.:
  {\em Mathematical Recreations and Essays}, 13th ed.
  New York: Dover, 1987.
\bibitem{Bol} Bol, G(errit?):
  Beantwoording van prijsvraag no.~17,
  {\em Nieuw Archief voor Wiskunde} {\bf 18} (1936), 14--66.
\bibitem{Brass} Brass, Peter:
  On point sets without $k$ collinear points, pages 185--192 in
  {\em Discrete Geometry: In honor of W.~Kuperberg's 60th Birthday}
  (A.~Bezdek, Ed.), Marcel Dekker Inc.,
  Pure and Applied Mathematics Series Vol.~{\bf 253}, 2003.
\bibitem{BMP} Brass, Peter; Moser, William; Pach, J\'anos:
  {\em Research Problems in Discrete Geometry.}
  New York: Springer, 2005.
\bibitem{Burr} Burr, Stefan A.:
  Planting Trees.  Pages 90--99 in
  {\em The Mathematical Gardner} (David Klarner, ed.;
  Belmont, Calif.: Wadsworth, 1981).
\bibitem{BGS} Burr, Stefan A.; Gr\"unbaum, Branko; Sloane, Neil J.A.:
  The Orchard Problem,
  {\em Geom.\ Dedicata} {\bf 2} (1974), 397--424.
\bibitem{Dudeney1} Dudeney, Henry E.:
  {\em The Canterbury Puzzles, and other curious problems.}
  New York: E.P.~Dutton \&~Co., 1908.
\bibitem{Dudeney2} Dudeney, Henry E.:
  {\em Amusements in Mathematics}.
  New York: Dover, 1958, 1970 (orig.\ Th.~Nelson \&~Sons, 1917).
\bibitem{Elkies} Elkies, Noam D.:
  ``New observations on an old puzzle:
  A plane configuration of $16$ points with $15$ lines of four''
  (November~2002),
  online at {\sf http://math.harvard.edu/$\sim$elkies/Misc/A5.html}\;.
\bibitem{EPS} Elkies, Noam D.; Pretorius, Lou M.; Swanepoel, Konrad J.:
  Sylvester-Gallai Theorems for Complex Numbers and Quaternions,
  to appear in {\em Discrete and Computational Geometry}
  ({\sf www.arxiv.org/math.MG/0403023}).
\bibitem{Grunbaum} Gr\"unbaum, Branko:
  New views of some old questions of combinatorial geometry,
  pages 451--468 in
  {\em Int.\ Teorie Combinatorie, Roma 1973, Tomo~I} (1976).
\bibitem{Hirzebruch} Hirzebruch, Friedrich E.P.:
  Arrangements of lines and algebraic surfaces.
  Pages 113--140 in {\em Arithmetic and Geometry, Vol. II}
  (= {\em Progr.\ Math.}\ {\bf 36}), Boston: Birkh\"auser, 1983.
\bibitem{Ismailescu} Ismailescu, Dan:
  Restricted point configurations with many collinear \hbox{$k$-tuplets},
  {\em Discrete Comput.\ Geom}\. {\bf 28} (2002), \#4 (571--575).
\bibitem{Kelly} Kelly, L.M.:
  A resolution of the Sylvester-Gallai problem of J.-P. Serre,
  {\em Discrete Comput.\ Geom.}\ {\bf 1} (1986), 101--104.
\bibitem{Palasti} Pal\'asti, Ilona:
  A construction for arrangements of lines
  with vertices of large multiplicity,
  {\em Studia Sci.\ Math.\ Hungar.}\ {\bf 21} (1986), \#1--2 (67--78).
\bibitem{PR} Poonen, Bjorn; Rubinstein, Michael: 
  The number of intersection points
  made by the diagonals of a regular polygon,
  {\em SIAM J. Discrete Math.}\ {\bf 11} (1998), \#1, 135--156
  ({\sf www.arxiv.org/math.MG/9508209}).
\bibitem{Vakil} Vakil, Ravi:
  {\em A Mathematical Mosaic: Patterns \& Problem Solving.}
  Burlington, Ontario: Brendan Kelly, 1997.
\end{thebibliography}
\end{document}